\documentclass[12pt,a4paper,twoside]{amsart}

\usepackage{mathrsfs}          
\usepackage[all]{xy}    
\usepackage[latin1]{inputenc}
\usepackage{amscd}
\usepackage{latexsym}
\usepackage{multicol}
\usepackage{amsmath}
\usepackage{amscd}
\usepackage{amssymb}
\usepackage{amsthm}

\theoremstyle{plain}
\newtheorem{Theorem}{Theorem}[section]
\newtheorem{lemma}[Theorem]{Lemma}
\newtheorem{corollary}[Theorem]{Corollary}
\newtheorem{proposition}[Theorem]{Proposition}

\newtheorem{conjecture2*}{Conjecture}
\newtheorem{remark2*}[conjecture2*]{Remark}
\newtheorem{Theorem2*}[conjecture2*]{Theorem}
\newtheorem{proposition2*}[conjecture2*]{Proposition}
\newtheorem{claim2*}[conjecture2*]{Claim}
\newtheorem{example2*}[conjecture2*]{Example}
\newtheorem{lemma2*}[conjecture2*]{Lemma}

\theoremstyle{definition}

\newtheorem{remark}[Theorem]{Remark}

\newtheorem{pkt}[Theorem]{}

\newtheorem{recall}[Theorem]{Recall}
\newcommand{\sB}{{\mathcal B}}
\newcommand{\sC}{{\mathcal C}}

\newcommand{\sE}{{\mathcal E}}
\newcommand{\sF}{{\mathcal F}}

\newcommand{\sM}{{\mathcal M}}

\newcommand{\sS}{{\mathcal S}}

\newcommand{\sW}{{\mathcal W}}

\newcommand{\sX}{{\mathcal X}}
\newcommand{\sY}{{\mathcal Y}}

\newcommand{\B}{{\mathbb B}}
\newcommand{\C}{{\mathbb C}}

\newcommand{\F}{{\mathbb F}}

\newcommand{\bP}{{\mathbb P}}
\newcommand{\Q}{{\mathbb Q}}
\newcommand{\R}{{\mathbb R}}

\newcommand{\id}{{\rm id}}

\newcommand{\Hg}{{\rm Hg}}

\newcommand{\SL}{{\rm SL}}
\newcommand{\GL}{{\rm GL}}

\title{Some Mirror partners with Complex multiplication}
\dedicatory{For Eckart Viehweg}
\author{Jan Christian Rohde}
\address{GRK1463/Institut f. Algebraische Geometrie\\ Leibniz Universit\"at Hannover\\ Welfengarten 1\\ 30167 Hannover\\ Germany}
\email{rohde@math.uni-hannover.de}
\begin{document}
\maketitle

\begin{abstract}
In this note we provide examples of families of Calabi-Yau 3-manifolds over Shimura varieties, whose mirror families contain subfamilies over Shimura varieties. Therefore these original families and subfamilies on the mirror side contain dense sets of complex multiplication fibers. In view of the work of S. Gukov and C. Vafa \cite{GV} this is of special interest in theoretical physics.
\end{abstract}

\section*{Introduction}
In theoretical physics rational conformal field theories are considered as particularly interesting class of conformal field theories. Let $(\sX,\sY)$ be a pair of families of Calabi-Yau 3-manifolds, which are mirror partners, $X$ be a fiber of $\sX$ and $Y$ be a fiber of $\sY$. In \cite{GV} S. Gukov and C. Vafa explain that $X$ and $Y$ yield a rational conformal field theory, if and only if both fibers have complex multiplication ($CM$). A family of Calabi-Yau manifolds over a Shimura variety has a dense set of $CM$ fibers, if the variation of Hodge structures ($VHS$) is related to the Shimura datum of the base space in a natural way as in \cite{JCR}. At present several of such families of Calabi-Yau 3-manifolds over Shimura varieties are known \cite{Bc}, \cite{G}, \cite{JCR}, \cite{JCR3}, \cite{VZ}. In general one does not know a Shimura subvariety of the base space on the mirror side.\footnote{Moreover for some of these examples \cite{G}, \cite{JCR3} the existence of a mirror is not clear.} Here we give new examples of pairs of families of Calabi-Yau 3-manifolds over Shimura varieties, which are subfamilies of mirror partners.

We start with a family $\sC_3$ of degree 3 covers of $\bP^1$ with 6 different ramification points over an open Shimura subvariety $\sM_3 \subset (\bP^1)^3$. By using the Fermat curve of degree 3 and $\sC_3$, one can construct a family of $K3$ surfaces with a non-symplectic involution over $\sM_3$ as described in \cite{JCR}, Section 8. The Borcea-Voisin construction yields a family $\sW$ of Calabi-Yau 3-manifolds, which has a dense set of $CM$ fibers. A. Garbagnati and B. van Geemen \cite{GG} have given a more general method to construct $K3$ surfaces, which yields the same $K3$ surfaces for the fibers of $\sC_3$. The latter method allows to construct $K3$ surfaces with non-symplectic involutions over the boundary of $\sM_3 \subset (\bP^1)^3$, where branch points of the fibers of $\sC_3$ collide. Here we show that the Borcea-Voisin construction yields a family of Calabi-Yau 3-manifolds over a Shimura subvariety contained in the boundary of the base space of $\sW$, whose fibers are its own Borcea-Voisin mirrors. Moreover here we find a Shimura surface on the boundary of the base space of $\sW$ such that the fibers of a family of Calabi-Yau 3-manifolds over this surface are Borcea-Voisin mirrors of the fibers of $\sW$. We will also see that these families contain dense sets of $CM$ fibers.

\section{Construction of $K3$ surfaces by automorphisms}

In this Section we recall the construction of $K3$ surfaces by the methods in \cite{GG}. For this construction we use the following families of curves:
\begin{enumerate}
 \item The family $\sC_1$ is the family of genus 2 curves given by
$$V(y^3-x_1(x_1-x_0)^2(x_1-\lambda x_0)^2x_0) \to \lambda \in \sM_1 := \bP^1\setminus \{0,1,\infty\}.$$
\item The family $\sC_2$ is the family of genus 3 curves given by
$$V(y^3-x_1(x_1-x_0)(x_1-\alpha x_0)(x_1-\beta x_0)x_0^2) \to (\alpha,\beta) \in \sM_2,$$
where
$$\sM_2 :=(\bP^1\setminus \{0,1,\infty\})^2 \setminus \{\alpha = \beta\}.$$
\item The family $\sC_3$ is the family of genus 4 curves given by
$$V(y^3-x_1(x_1-x_0)(x_1-\alpha x_0)(x_1-\beta x_0)(x_1-\gamma x_0)x_0) \to (\alpha,\beta,\gamma) \in \sM_3,$$
where
$$\sM_3 :=(\bP^1\setminus \{0,1,\infty\})^3 \setminus (\{\alpha = \beta\}\cup\{\alpha = \gamma\}\cup\{\beta = \gamma\}).$$
\end{enumerate}

\begin{remark} \label{1.1}
The families $\sC_1$ and $\sC_2$ can be obtained by collision of the branch points of the fibers of $\sC_3$ over the boundary divisor of $\sM_3 \subset (\bP^1)^3$.
Let $\Gamma$ denote the monodromy group of the $VHS$ of $\sC_3$. Note that for $j = 1,2,3$ one can apply the Deligne-Mostow theory \cite{DM} to the $VHS$ of $\sC_j$. For an overview of this topic see also \cite{Loo}. Due to the Deligne-Mostow theory, the period domain of the family $\sC_j$ is the complex ball $\B_j$ and $\sM_j$ is a dense open subset of $\Gamma\backslash \B_j$. In this sense the base spaces $\sM_j$ are modular. Moreover $\sM_2$ and $\sM_1$ are contained in the complement of $\sM_3$ in $\Gamma\backslash \B_3$ (follows from \cite{Loo}, Theorem $3.1$ and the description of the period map in \cite{Loo}, Section 4).

One can also see that $\sM_j$ is an open dense subset of a Shimura variety, which is a ball quotient. This can be concluded from the type of $VHS$ of the given families (compare \cite{JCR}, Subsection $6.3$) and the description of such a $VHS$ in the proof of \cite{JCR}, Theorem $4.4.4$ in combination with the description of the period map above.
\end{remark}

For $j = 1,2,3$ and $p \in \sM_j$ let $f_j(t)\in \C[t]$ be a degree 6 polynomial such that $(\sC_j)_p$ is given by the equation $v^3-f_j(t)=0$. Moreover let $\xi = e^{2\pi i\frac{1}{3}}$. It is clear that $\sC_j$ has the $\sM_j$-automorphism fiberwise given by
$$\beta_j:(v,t) \to (\xi v, t).$$
Let
$$\F_3 = V(y^2z-x^3-z^3)\subset \bP^2$$
be a genus 1 curve isomorphic to the Fermat curve of degree 3 and
$$\alpha_{\F_3}:\F^3 \to \F^3 \ \ \mbox{be given by} \ \ (x:y:z) \to (\xi x:y:z).$$
We have chosen this explicite formula due to technical reasons. Moreover let $S_{f_j}$ be a minimal model of a surface given by the Weierstrass equation
$$Y^2 = X^3 +f_j^2(t).$$
For the following lemma we will use methods, which occur already in the proof of \cite{GG}, Proposition $2.2$:

\begin{lemma}
The surface $S_{f_j}$ is a $K3$ surface birationally equivalent to $\F_3 \times (\sC_j)_p/(\alpha_{\F_3},\beta_j)$.
\end{lemma}
\begin{proof}
The rational map
$$m_j:\F_3 \times (\sC_j)_p \to S_{f_j}$$
is given by
$$((t,v),(x,y)) \to (v^2x,v^3y,t).$$
The reader checks easily that $m_j$ is $(\alpha_{\F_3},\beta_j)$-invariant and of degree 3. Moreover one computes
$$Y^2 = (v^3y)^2 = v^6y^2=v^6(x^3+1)= (v^2x)^3+ f_j^2(t) = X^3+f_j^2(t).$$
From \cite{GG} we know that the minimal model $S_{f_j}$ is a $K3$ surface.
\end{proof}

\section{Some Automorphisms of our $K3$ surfaces}

\begin{pkt}
The surface $S_{f_j}$ has an elliptic fibration given by
$$S_{f_j} \to \bP^1 \ \ \mbox{via} \ \ (X,Y,t) \to t$$
in the following way (see also \cite{GG}):

If $f_j(t_0) \neq 0$, the fiber of $t_0$ is given by the elliptic curve
$$V(Y^2Z = X^3 +uZ^3) \subset \bP^2, \ \ \mbox{where} \ \ u = f_j^2(t_0).$$

Now let $t_0 \in \bP^1$ be a zero of $f_j(t)$. By using the Tate algorithm, one can compute the singular fibers. If $f_j(t)$ has a simple zero in $t_0$, the singular fiber $(S_{f_j})_{t_0}$ is of type ${\bf IV}$. Thus it consists of three rational curves intersecting transversally in one point.

Now assume that $f_j(t)$ has a double zero in $t_0$. Then the fiber $(S_{f_j})_{t_0}$ is of type ${\bf IV}^*$. Thus it is given by 7 rational curves with the following intersection graph of type $\tilde E_6$:
$$
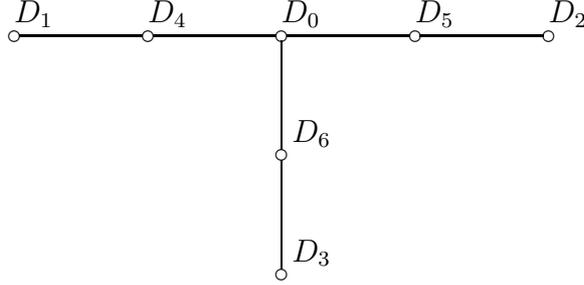
\begin{picture}(200, 100)

\put(2,90){\line(1,0){46}}
\put(52,90){\line(1,0){46}}
\put(102,90){\line(1,0){46}}
\put(152,90){\line(1,0){46}}

\put(100,88){\line(0,-1){41}}
\put(100,43){\line(0,-1){41}}

\put(100,95){$D_0$}
\put(100,90){\circle{4}}

\put(0,95){$D_1$}
\put(0,90){\circle{4}}

\put(200,95){$D_2$}
\put(200,90){\circle{4}}

\put(104,5){$D_3$}
\put(100,0){\circle{4}}

\put(50,95){$D_4$}
\put(50,90){\circle{4}}

\put(150,95){$D_5$}
\put(150,90){\circle{4}}

\put(104,50){$D_6$}
\put(100,45){\circle{4}}

\end{picture}
$$
\end{pkt}

Let $\iota_{f_j}:S_{f_j}\to S_{f_j}$ denote the involution given by 
$$(X,Y,t) \to (X,-Y,t)$$
and $\alpha_{f_j}:S_{f_j}\to S_{f_j}$ denote the automorphism of degree 3 given by
$$(X,Y,t) \to (\xi X,Y,t).$$

\begin{pkt} \label{invy}
The fixed locus of $\iota_{f_j}$ contains clearly the section $s_{\infty}$ of the elliptic fibration fiberwise given by $(0:1:0) \in \bP^2$ for a general $t \in \bP^1$ and the curve $C_{j,p}$ given by
$$X^3+f_j^2(t) = 0,$$
which is isomorphic to $(\sC_j)_p$ (see \cite{JCR}, Remark $2.1.8$). One can also verify by explicit computation that one obtains a fiber $F$ of type ${\bf IV}$ over a simple zero of $f_j(t)$ by blowing up once. This computation shows that $\iota_{f_j}$ interchanges two irreducible components of $F$ and $C_{j,p}$ intersects $F$ in the intersection point of its irreducible components. The involution $\iota_{f_j}$ acts non-trivially on the third irreducible component of $F$.

Moreover the fixed locus of the automorphism $\alpha_{f_j}$ contains also the section $s_{\infty}$ and the sections $s_{\pm}(t) = (0:\pm f_j(t):1)\in \bP^2$.

Since the fixed loci of $\iota_{f_j}$ and $\alpha_{f_j}$ contain curves, $\iota_{f_j}$ and $\alpha_{f_j}$ are non-symplectic.
\end{pkt}

\begin{recall}
It is well-known that a non-symplectic involution of a $K3$ surface has a fixed locus, which is either empty or consists of smooth disjoint curves. Moreover it is well-known that a non-symplectic automorphism of degree 3 of a $K3$ surface has a fixed locus, which is empty or consists of smooth disjoint curves and finitely many isolated points. We will use these facts later.
\end{recall}

It is a little bit harder to compute the intersection of the fixed locus of $\iota_{f_j}$ with a singular fiber $F^*$ of type ${\bf IV^*}$. For doing this we also consider the automorphism $\alpha_{f_j}$.

First we can a priori state that either $\iota_{f_j}|_{D_0} = \id$ or that one has without loss of generality $\iota_{f_j}(D_0\cap D_4) = D_0\cap D_5$. In both cases one obtains without loss of generality that $\iota_{f_j}(D_3 \cup D_6) = D_3 \cup D_6$. Hence on $D_3 \cup D_6$ one finds without loss of generality an isolated fixed point of $\iota_{f_j}|_{F^*}$, which is an intersection point with the section $s_{\infty}$ or $C_{j,p}$. This point is also fixed by $\alpha_{f_j}$. Thus the automorphism $\alpha_{f_j}$ has to act as a permutation $\sigma\in A_3$ on the intersection points $D_0 \cap D_k$ for $k = 4,5,6$, which fixes $D_0\cap D_6$. Hence $\sigma = \id$. The Hurwitz formula tells us that either $\alpha_{f_j}|_{D_0}= \id$ or that the quotient map by the degree 3 automorphism $\alpha_{f_j}|_{D_0}$ has 2 ramification points. Thus $\alpha_{f_j}|_{D_0}= \id$. Hence for $k = 4,5,6$ the fixed points of $\alpha_{f_j}|_{D_k}$ are given by the intersection points of $D_k$ with the other irreducible components of the fiber $F^*$.

Since $s_{\infty}$ is fixed by $\alpha_{f_j}$, we assume without loss of generality that $s_{\infty}$ hits a singular fiber $F^*$ of type ${\bf IV}^*$ in $D_3$. 
%Thus $\alpha_{f_j}(D_3) = D_3$ and one fixed point of this automorphism on $D_3$ is given by $D_3\cap D_6$.
Moreover one has also that $\iota_{f_j}(D_3) = D_3$ and $\iota_{f_j}$ acts as a non-trivial involution on $D_3$ with a fixed point $s_{\infty}\cap D_3$.

\begin{lemma} \label{star}
One cannot have $\iota_{f_j}|_{D_0} = \id$.
\end{lemma}
\begin{proof}
Assume that $\iota_{f_j}|_{D_0} = \id$. Hence $\iota_{f_j}|_{D_6} \neq \id$. Since $s_{\infty}$ hits $F^*$ in $D_3$, one has a non-trivial involution on $D_3$, whose fixed points are given by the intersection points with the fixed curves $s_{\infty}$ and $C_{j,p}$. Thus $\iota_F$ acts on $D_1$ and $D_2$ as the identity map. Since the sections $s_{\pm}$ are interchanged by $\iota_{f_j}$, one concludes $s_{\pm}\cap F^* \notin D_0, D_1, D_2$. Since $\alpha_{f_j}$ fixes the points $D_0\cup D_4$ and $D_1\cup D_4$ respectively and acts non-trivially on $D_4$, the point $s_{\pm}\cap F^*$, which is fixed by $\alpha_{f_j}$, cannot be contained in $D_4$. By analogue arguments, one concludes $s_{\pm}\cap F^* \notin D_5, D_6$. Since the section $s_{\infty}$, which is also fixed by $\alpha_{f_j}$, hits $D_3$, it is not possible that both sections $s_{\pm}$ intersect $D_3$, on which $\alpha_{f_j}$ acts non-trivially. Contradiction!
\end{proof}

The involution $\iota_{f_j}|_{F^*}$ has an isolated fixed point $s_{\infty}\cap D_3$. Since the intersection point $F^*\cap C_{j,p}$ is the only additional isolated fixed point of $\iota_{f_j}|_{F^*}$, one cannot have 3 isolated fixed points on $D_3\cup D_6$ with respect to $\iota_{f_j}$. Thus one concludes:

\begin{corollary}
The curve $D_6$ is contained in the fixed locus with respect to $\iota_{f_j}$ and $\iota_{f_j}$ interchanges the handle consisting of $D_1$ and $D_4$ with the handle consisting of $D_2$ and $D_5$.
\end{corollary}

\section{Construction of mirror pairs with complex multiplication}

\begin{recall}
Let $S$ be a $K3$ surface with non-symplectic involution $\iota$, which has a fixed locus consisting of the curves $C_1, \ldots, C_N$, and $E$ be an elliptic curve with involution $\iota_E$ fixing 4 points. Moreover let
$$N^{\prime} = \sum\limits_{i=1}^N g(C_i),$$
where $g(C_i)$ denotes the genus of $C_i$. Then the Calabi-Yau 3-manifold $X$ obtained from the Borcea-Voisin construction given by blowing up the singularities of $S\times E/(\iota_S,\iota_E)$ once has the Hodge numbers
\begin{equation} \label{vo}
h^{1,1}(X) = 11 + 5N - N^{\prime} \ \ \mbox{and} \ \ h^{2,1}(X) = 11 + 5N^{\prime} - N
\end{equation}
(see \cite{Voi2}).

In many cases the involution on $H^2(S,X)$, which is given by the action of $\iota$, can be used to construct a second involution $\iota'$ on the lattice $H^2(S,X)$. The involution $\iota'$ can be realized as an involution of a family of $K3$ surfaces $\sS'\to \sB$ over $\sB$, whose restrictions to each fiber of $\sS'$ are non-symplectic involutions. By a relative version of the construction above for $\sS'$, one obtains the Borcea-Voisin mirror family of $X$ (for details see \cite{Bc2}, \cite{Voi2}).
\end{recall}

\begin{remark}
By using $\sC_3$, one has already constructed families of Calabi-Yau manifolds over Shimura varieties (see \cite{JCR}, Section 8 and \cite{JCR}, Section 9). For this construction in \cite{JCR} one has used a family of $K3$ surfaces, which is precisely the family over $\sM_3$, which occurs in \cite{GG}, Remark $1.3$ and also here. This follows from the fact that both constructions yield a $K3$ surface, which is a minimal model of $\F_3 \times (\sC_j)_p/(\alpha_{\F_3},\beta_j)$. Despite this fact it is not clear that the family $\sX_3$ of Calabi-Yau 3-manifolds, which we construct below, is the family $\sW$ in the notation of \cite{JCR}. Nevertheless both families are contained in the same Borcea-Voisin family. Thus the fibers have the same Borcea-Voisin mirrors.
\end{remark}

\begin{pkt}
By $\ref{invy}$, the involution $\iota_{f_j}$ on $S_{f_j}$ has a fixed locus containing a rational curve $s_{\infty}$ and the curve $C_{j,p}$ of genus $j+1$. Moreover the elliptic fibration of $S_{f_j}$ contains $3-j$ singular fibers of type ${\bf IV}^*$ and each of these fibers has one rational curve contained in the fixed locus of $\iota_{f_j}$ (see Corollary $\ref{star}$). Thus by using $\eqref{vo}$ and the family of elliptic curves
$$\sE \to \sM_1, \ \ V(y^2z-x(x-z)(x-\lambda z)) \to \lambda,$$ the Borcea-Voisin construction yields families $\sX_j \to \sM_j\times \sM_1$ of Calabi-Yau 3-manifolds with the following Hodge numbers:
$$\begin{tabular}{|c||c|c|} \hline
$j$ & $h^{1,1}$ & $h^{2,1}$ \\ \hline \hline
$3$ & 17 & 29 \\ \hline
$2$ & 23 & 23 \\ \hline
$1$ & 29 & 17\\ \hline
\end{tabular}
$$
\end{pkt}

\begin{remark}
By \cite{Bc2}, Section 3 and Section 4, one can easily check that $\sX_2$ is contained in a family, which is its own Borcea-Voisin mirror family. Moreover the families $\sX_1$ and $\sX_3$ can be embedded in families, which are Borcea-Voisin mirrors of each other.
\end{remark}

\begin{remark} \label{3.3}
By the construction above, the families $\sX_1$ and $\sX_2$ are contained in the boundary of $\sX_3$. Moreover by using Remark $\ref{1.1}$,
one can show that the period map of $\sX_j$ is a multivalued map to a dense open subset of $\B_j\times \B_1$. From these results one can conclude that the base space of $\sX_j$ is an open subset of a Shimura variety with associated Hermitian symmetric domain $\B_j\times \B_1$.

By analogue arguments, one can also see that $\sX_1$ is defined over the boundary of $\sX_2$.
\end{remark}

By \cite{Voi2}, $2.21$, we have a precise description how a fiber of $\sX_j$ provides a $(1,1)$-form on a fiber of $\sX_{4-j}$ by the mirror map. Due to \cite{GV} one can assume that each pair of complex multiplication fibers of $\sX_{j}$ and $\sX_{4-j}$ yields a rational conformal field theory. Now we are going to show that each $\sX_j$ has a dense set of complex multiplication ($CM$) fibers for $j = 1,2,3$. First recall the definition of $CM$:

\begin{pkt}
Let $X$ be a compact K\"ahler manifold of complex dimension $n$ and $S^1$ be the $\R$-algebraic group
$$S^1 = {\rm Spec}(\R[x,y]/x^2+y^2 - 1),$$
where
$$S^1(\R)=\left\{M = \left(\begin{array}{cc}
a & b \\
-b  & a \end{array}\right) \in \SL_2(\R) \right\}\cong \{z \in \C:|z| = 1\}.$$
The rational Hodge structure on $H^n(X,\Q)$ of weight $n$ corresponds to the representation
$$h_X:S^1\to \GL(H^n(X,\R)), \ \ h_X(z)v =z^p\bar z^qv \ \  (\forall v \in H^{p,q}(X) \ \ \mbox{with} \ \ p + q = n).$$
The Hodge group $\Hg(X)$ is the smallest $\Q$-algebraic subgroup $G$ of $\GL(H^n(X,\Q))$ such that $h_X(S^1) \subset G_{\R}$.
We say that $X$ has $CM$, if $\Hg(X)$ is a torus.

For more details in the case of Calabi-Yau 3-manifolds see \cite{Bc}.
\end{pkt}

\begin{proposition}
For $j = 1,2,3$, the family $\sX_j$ has a dense set of $CM$ fibers.
\end{proposition}
\begin{proof}
By \cite{JCR}, Subsection $6.3$, each family $\sC_j$ has a dense set of $CM$ fibers. Note that the family of elliptic curves
$$\sE \to \sM_1, \ \ V(y^2z-x(x-z)(x-\lambda z)) \to \lambda$$
has also a dense set of $CM$ fibers. Since the ramification locus of the involution $\iota_{f_j}$ on $S_{f_j}$ consists of $C_{j,p}\cong (\sC_j)_p$ and some rational curves, it remains to show that $S_{f_j}$ has $CM$, if $(\sC_j)_p$ has $CM$. Using this result one can then conclude as in \cite{JCR}, Subsection $7.2$ that $(\sX_j)_{(p,q)}$ has $CM$, if $(\sC_j)_p$ and $\sE_q$ have $CM$.

The singularities of the fibers of $\F_3 \times \sC_j/(\alpha_{\F_3},\beta_j)$ are given by the singular sections. Let $m_j$ denote the quotient map by $(\alpha_{\F_3},\beta_j)$. Near the sections of fixed points corresponding to the singular sections of $\F_3 \times \sC_j/(\alpha_{\F_3},\beta_j)$ the action of $(\alpha_{\F_3},\beta_j)$ is given by $(\xi, \xi)$ or $(\xi, \bar \xi)$.

First consider the case $(\xi, \bar \xi)$. In this case one blows up the corresponding sections on $\F_3 \times \sC_j$ with exceptional divisor $E_1$. The automorphism $(\alpha_{\F_3},\beta_j)$ does not act trivially on $E_1$. Thus we blow up the two fixed sections on each connected component of $E_1$ with smooth exceptional divisor $E_2$. This divisor is contained in ramification locus of $m_j$. Now the quotient by $(\alpha_{\F_3},\beta_j)$ is smooth in a neighbourhood of $m_j(E_1\cup E_2)$.

In the case $(\xi, \xi)$ we blow up the section of fixed points and obtain a smooth exceptional divisor contained in the ramification locus. 

Let $\widetilde{\F_3\times \sC_j}$ denote the manifold obtained from the previous blowing up operations on $\sC_j\times \F_3$ and 
$$\sF_j = \widetilde{\F_3\times \sC_j}/(\alpha_{\F_3},\beta_j).$$
Thus we obtain a model $\sF_j$ of the quotient $\F_3\times (\sC_j)_p/(\alpha_{\F_3},\beta_j)$ consisting of smooth fibers over $\sM_j$. The surface $\F_3\times \sC_j$ has $CM$, if $\F_3$ and $(\sC_j)_p$ have $CM$. Note that monoidal transformations of surfaces do not have any effect to the property of $CM$ (compare \cite{JCR}, Corollary $7.1.6$). Since the Hodge structure on $H^2((\sF_j)_p,\Q)$ is a sub-Hodge structure of the Hodge structure on $H^2(\widetilde{(\F_3\times \sC_j)}_p,\Q)$, one concludes that $(\sF_j)_p$ has $CM$, if $(\sC_j)_p$ has $CM$. Moreover we can use monoidal transformations to obtain $S_{f_j}$ from $(\sF_j)_p$. Thus $S_{f_j}$ has $CM$, if $(\sF_j)_p$ has $CM$.
\end{proof}

\section*{Acknowledgements}
This paper was written at the Graduiertenkolleg ``Analysis, Geometry and String Theory'' at Leibniz Universit\"at Hannover. I would like to thank Lars Halle for instructive discussions about Neron models and for once asking the question ``What happens on the boundary of $\sW$?''. This question is partially answered here. I would like to thank Klaus Hulek and Matthias Sch\"utt for their interest and their comments, which helped to improve this text. Moreover I would like to thank Bernd Siebert, who has pointed out that it is of interest to understand the period maps of the given families, which are outlined in Remark $\ref{1.1}$ and Remark $\ref{3.3}$.

\end{document}